\documentclass[12pt,a4paper,twoside,final,notitlepage, leqno]{article}
\usepackage[english]{babel}
\usepackage[T1]{fontenc}  
\usepackage{epsfig, graphicx, amssymb}
\usepackage{amsmath,amsthm,epsfig,amsfonts,bbm}  
\usepackage[bottom]{footmisc}
 
\def\br {\break}

\newcommand{\moneq}{\vspace*{-6pt} \begin{equation} \displaystyle } 
\newcommand{\moneqstar}{\vspace*{-6pt} \begin{equation*} \displaystyle } 
\newcommand{\monendstar}{\vspace*{-6pt} \end{equation*}   }
\newcommand{\monend}{\vspace*{-6pt} \end{equation}   }

\def\section*#1{}

\usepackage{fancyhdr}
\renewcommand{\headrulewidth}{0pt}

\begin{document} 

\fancypagestyle{plain}{ \fancyfoot{} \renewcommand{\footrulewidth}{0pt}}
\fancypagestyle{plain}{ \fancyhead{} \renewcommand{\headrulewidth}{0pt}}

~

\bigskip \bigskip   \bigskip  
 
\centerline {\bf \LARGE Kron's method and cell complexes } 

 \bigskip  \centerline {\bf \LARGE for magnetomotive and electromotive forces }

\bigskip \bigskip 

\centerline { \large  Olivier Maurice$^{a}$ , Alain Reineix$^{b}$, 
Philippe Durand$^{c}$, Fran\c{c}ois Dubois $^{d,e}$ }   

\bigskip  

\centerline { \it  \small   
$^a$     GERAC,  3 avenue Jean d'Alembert, ZAC de Pissaloup, 78190 Trappes, France.}

\centerline { \it  \small   
$^b$     Limoges University, XLIM - UMR CNRS n$^{\rm o}$7252, }

\centerline { \it  \small 123, avenue Albert Thomas, 87060 Limoges, France.}

\centerline { \it  \small   
$^c$   Conservatoire National des Arts et M\'etiers, Department of Mathematics, } 

\centerline { \it  \small 
Mod\'elisation Math\'ematique et Num\'erique, 292 rue Saint Martin, 75141 Paris, France.}

\centerline { \it  \small   
$^d$  Conservatoire National des Arts et M\'etiers, Structural Mechanics and Coupled Systems } 

\centerline { \it  \small 
 Laboratory, 292 rue Saint Martin, 75141 Paris, France} 

\centerline { \it  \small   
$^e$  University  Paris-Sud, Department of mathematics, Orsay, France} 


\smallskip

\centerline { \it  \small 
olivier.maurice@gerac.com,  alain.reineix@xlim.fr, }

\centerline { \it  \small 
   philippe.durand@cnam.fr, francois.dubois@math.u-psud.fr }

\bigskip  

\centerline { {\rm  28 november  2014} 
\footnote {\rm  \small $\,\,$ Article published in 
{\it IAENG International Journal of Applied Mathematics}, volume~44, 
issue~4, pages 183-191, november 2014. }} 

\bigskip  

\noindent  {\bf Abstract  } 

\noindent   
Starting from topological principles we first recall the elementary ones giving Kirchhoff's 
laws for current conservation. Using in a second step the properties of spaning tree, 
we show that currents are under one hypothesis intrinsically boundaries of surfaces flux. 
Naturally flux appears as the object from which the edge comes from. The current becomes the 
magnetomotive force (mmf) that creates the flux in the magnetostatic representation. Using a 
metric and an Hodge's operator, this flux creates an electromotive force (emf). This emf is 
finally linked with the current to give the fundamental tensor - or ``metric'' - of the Kron's 
tensorial analysis of networks. As it results in a link between currents of cycles 
(surface boundaries) and energy sources in the network, we propose to symbolize this 
cross talk using chords between cycles in the graph structure on which the topology is based. 
Starting then from energies relations we show that this metric is the Lagrange's operator 
of the circuit. But introducing moment space, the previous results can be extended to non 
local interactions as far field one. And to conclude, we use the same principle to 
create general relation of information exchange between networks as functors between categories.



\newpage

\fancyhead[EC]{\sc{Olivier Maurice, Alain Reineix, Philippe Durand, Fran\c{c}ois Dubois}} 
\fancyhead[OC]{\sc{Kron's method and cell complexes\ldots }} 
\fancyfoot[C]{\oldstylenums{\thepage}}

\bigskip \bigskip  \noindent  {\bf \large 1) \quad  Introduction }   

\noindent  
Gabriel Kron (1900 - 1968) has transfered the tensorial analysis developed 
in the framework of general relativity to the world of applied electromagnetism  
\cite{GKRONI,GUP}. He had felt the straight relations between Kirchhoff's laws and topology. 
To replace his reflexions at this time, we try to give a brief history of these works, 
from Kron to nowdays. 

\noindent  
As many words exist to call the various elements available on graphs, 
we propose some of them that we use in our paper. Our first purpose is to find again 
these straight relations as simplest as possible. Making this exercise and under the hypothesis 
that no spaning tree current sources are involved, only electrostatic and magnetostatic 
phenomenons are considered, it appears that the current can be seen as a fundamental 
element of the faces boundary space. Through this vision, the electromotive force belongs 
to the faces space. The current seen as a boundary of faces, becomes equivalent to 
magnetomotive forces. By a physical understanding, a relation between magnetomotive 
forces and electromotive forces must exists. A problem appears because both doesn't 
belong to the same differential form dimension. Thanks to the Hodge operator, 
we make this link, following previous existing works. But if the electromotive force, 
which give the network its energy, belongs to the faces space, self inductance reaction 
must belongs to the same space and mutual inductance interaction translates the 
relation of Hodge. 

\noindent  
These components, as sensed by Kron, represent the metric that 
we take to compute these interactions. This leads to the Lagrangian expression of 
the whole graph and to the ``chords'' elements that symbolize these interactions. 
This Lagrangian must then be increased taking into account the spaning tree sources, 
added to the faces ones. This new space give us the complete base to take into account 
far field interactions. Lamellar fields create current sources, rotationnal fields 
create electromotive forces. This continuous fields are connected to our topology 
using moment space. This space is the frontier between our first bounded manifolds 
which generate the graph and continuous not bounded manifolds which are radiated and 
propagated fields. The first topological discussion gives all the base to include 
this new interaction through a generalized definition of the chords. At each step, 
we first give a topological approach before to ``translate'' it in expressions 
more usually given by physicists.

\bigskip \bigskip  \noindent  {\bf \large 2) \quad  
Tensorial Analysis of Networks (TAN) history }   

\noindent 
Gabriel Kron has written is famous ``tensorial analysis of networks'' in 1939. 
Before this work, he has written in 1931 a first remarkable paper ``Non-Riemannian 
Dynamics of Rotating Electrical Machinery''  \cite{K39}. For this work he had the 
Montefiore  price of the university of Liege in 1933 and the M.I.T. journal of 
mathematics and physics publishes the entire paper in the May 1934  \cite{K34}. 
This paper instantly produces wide-spread discussion and contreversy \cite{GKAST}. 
Kron uses its own notation without regarding established ones. This leads to some 
mathematicians contempt. But some of them were clear enought to understand and study 
Kron's work. Hoffmann \cite{BAH}, Roth \cite{ROTI} make links between Kron's 
concept and topology ones. Physicists like Branin, Happ, and in France Denis-Papin and 
Kaufmann \cite{PEK} promote Kron's work for electrical engineers. Many studies 
were done after around the concept of \textit{Diakoptic} initiated by Kron \cite{HAPP}. 
But to focus on topology, less references are available. Major lecture was made by 
Balasubramanian, Lynn and Sen Gupta \cite{GUP}. Recently, there is the work done by 
Gross and Kotiuga \cite{GK}, following first one of Bossavit \cite{BOS}. These last 
two works was made more particularly for finite element method. But they give 
fundamental bases through algebraic topology, following previous works of Roth and 
others, clarified using benefit of years passing. In this paper we try to take 
benefit of all this story and to present as clear as we can our understanding 
on these concept. 

\bigskip \bigskip  \noindent  {\bf \large 3) \quad  Notations}

\noindent 
We note $\mathbb{N}$ the set of integers, $\mathbb{R}$ the set of real numbers, 
$\mathbb{C}$ for complex, $\mathcal{T}$ for cells, etc. We work in a complex 
cellular $\mathcal{T}^{\infty}$ made of vertexes $s\in \mathcal{T}^0$, edges 
$a\in \mathcal{T}^1$ and faces $f\in \mathcal{T}^2$, etc. $\mathcal{T}$ is 
the whole set of these geometrical or chain objects. Low indices refer to chains, 
high ones to the geometric objects. Geometric objects are classical forms. 
Chains are abstract objects embedding properties added to the geometric objects 
in order to represent symbolically a real thing. A set of currents runing in a 
system can be linked with a set of chains associated with edges. This set of currents 
are components of a unique current vector. The current vector constitute a chain, 
image of some real currents on an electronic system. More generally, vectors can 
be associated with each of these geometric objects and their bases: 
$\left| s\right\rangle$ for vertexes, $\left| a\right\rangle$ for edges, $\left|f\right\rangle$ 
for faces, {\it etc}.
We note $\mathcal{T}_j$ the vectorial space created by the geometric vectors 
$\left|x\right\rangle$ of $\mathcal{T}^{D(x)}$ ($D(x)=0$ if x similar to s, {\it etc.}). 
$\mathcal{T}_j$ can be developped as:
\begin{equation}
\mathcal{T}_j=\left\lbrace\sum_{\sigma \in \mathcal{T}^j} \alpha_\sigma 
\left|\sigma\right\rangle,~\alpha_\sigma\in\mathbb{R},~\mathbb{C}\right\rbrace,~j\in N
\end{equation}
In this definition, we see here a generalized formulation of the classical writing 
of a vector, using the mute index notation (each time an index is repeated, the 
summation symbol on the index can be omitted) in  \cite{PEIK}: 
$\vec{f}=f^\sigma {\vec{u}}_\sigma$. Here $\vec{f}$ is a vector developed on the 
base ${\vec{u}}_\sigma$ of components $f^\sigma$.

\bigskip \bigskip  \noindent  {\bf \large 4) \quad  Boundary operator}

\noindent 
We now introduce the boundary operator. It translates the intuitive understanding 
of object boundary. The boundary of a segment is a pair of two points, the one of 
a surface is a closed line, and so on. The boundary operator is the base of all 
Whitney's concepts  \cite{WHIT}. To define integration through bounded objects, 
anyone needs boundaries. Once more, this boundary concept is natural. It is an 
application from $\mathcal{T}$ to $\mathcal{T}$, $\partial :\mathcal{T}\rightarrow \mathcal{T}$, 
and more precisely, $\partial$ is an operator from ${\mathcal{T}}_j$ to ${\mathcal{T}}_{j-1}$. 
Its self composition leads to zero: 
$\forall\theta\in \mathcal{T}^{D(\theta)},~\partial\circ\partial \left|\theta\right\rangle=0$. 
For example, $\partial\circ\partial \left|a\right\rangle=0$: the boundary of an edge 
is a vertex and the boundary of a vertex is null, or 
$\partial\circ\partial \left|f\right\rangle=0$: the boundary of a face is a closed line, 
and the boundary of this closed line is null (remember $\partial\circ\partial=0$). 
As we will see in next paragraph, this operator is linked to various connectivities 
in tensorial algebra. A face boundary is a cycle (a closed line). When we have an 
edge $\left|a\right\rangle$, we can consider its boundary $\partial\left|a\right\rangle$ 
which is a couple of vertices. The boundaries can be developed on the zero chain 
of vertices :
\begin{equation}\label{Bord-Arete}
\partial\left|a\right\rangle=\sum_{s\in{\mathcal{T}}^0} B^s_{a} \left|a\right\rangle \,.
\end{equation}
Take a look to the graph figure 1. We can easily find its incidence B making 
relations between the vertices $s$ and the edges $a$:
\begin{equation}
B=\left[
\begin{array}{ccccc}
1&-1&-1&0&0\\-1&1&0&1&1\\0&0&1&-1&-1
\end{array}\right] \,.
\end{equation}
Each row is linked to an vertex and each column to an edge.

\bigskip \bigskip  \noindent  {\bf \large 5) \quad  Seeing electrical current as a 1-chain 
and the potential}  

\noindent \qquad \quad  {\bf \large    as a 0-cochain}

\noindent 
In the following we consider the electrical current $i$ as an element $\left|i\right\rangle$ 
of the space ${\mathcal{T}}_1$: on each edge $k$, the current has a component $i^k$ which 
is a real number:
\begin{equation}
\left|i\right\rangle= \sum_{k\in{\mathcal{T}}^1} i^k\left|k\right\rangle \,.
\end{equation}
In the tensorial analysis of networks \cite{GKRONI} (as previously in classical nodal 
techniques \cite{PEIK}), the boundary operator applied to edges is called the ``incidence". 
Using the sign rule saying that a current entering a vertex is affected of a plus sign and 
a current leaving a vertex is affected of a minus sign, it is a matrix that gives the 
relations between vertices and edges.
It is possible to create for each vertex $s$ a linear form $\left\langle s\right|$ acting 
on all the vertices: $\left\langle s|\sigma\right\rangle=0$ if $s$ and $\sigma$ differ, 
$\left\langle s|\sigma \right\rangle=1$ if $s=\sigma$.
This form belongs to ${\mathcal{T}}^{*}_0$: the dual space of ${\mathcal{T}}_0$ composed 
by the 0-cochains.
Moreover to each vertex $s\in{\mathcal{T}}^0$ is associated a potential value $V_s$. 
With this set of numbers we construct the potential $V$ such that:
\begin{equation}
\left\langle V\right|=\sum_{s\in{\mathcal{T}}^0} V_s  \left\langle s\right| \,.
\end{equation}

\bigskip \bigskip  \noindent  {\bf \large 6) \quad  A topological form of Kirchhoff's laws }

\noindent 
We propose here to write the Kirchhoff's laws in a single abstract form as:
\begin{equation}\label{VDI}
\left\langle V | \partial i\right\rangle=0,~\forall 
V\in{\mathcal{T}}_0^{*},~\forall i\in{\mathcal{T}}^1 \,.
\end{equation}
This relation applied to electricity is usual for physicists. If we retain one 
particular node (for example node 1 on figure 1), the algebraic sum of the 
currents $i^k$ for $k$ an edge that contains the vertex number 1 is equal to zero; 
in this case $i^1-i^2+i^3=0$.
To interpret relation (\ref{VDI}) as the second Kirchhoff's law relative to the mesh law, 
we need the mathematical notion of co-boundary.

\noindent 
From the chains ${\mathcal{T}}_j$, we introduce the space ${\mathcal{T}}^{*}_j$ of 
co-chains of degree $j$: we have defined the duality product 
$\left\langle s | \sigma \right\rangle$ for two vertices $s$ and $\sigma$ of the 
cellular complex. We do the same for each geometrical object of dimension $j$. 
For $\alpha \in {\mathcal{T}}^j$, the dual form $\left\langle \alpha \right|$ belongs 
to ${\mathcal{T}}^{*}_j$ and is defined for each  $\theta \in {\mathcal{T}}^j$ by 
$\left\langle \alpha|\theta\right\rangle=0$ if $\alpha$ and $\theta$ differ, 
$\left\langle \alpha|\theta \right\rangle=1$ if $\alpha=\theta$. The co-boundary 
operator $\partial^{\rm o}$ is the polar operator of the boundary operator $\partial$. 
By definition for $ \varphi \in   {\mathcal{T}}_j^* $ and $\theta \in   {\mathcal{T}}_{j+1}$ 
we have:
\begin{equation}
\left\langle\partial^{\rm o} \varphi | \theta \right\rangle \equiv \left\langle  
\varphi |  \partial \theta\right\rangle, \quad \forall \varphi\in {\mathcal{T}}^{*}_j, 
\, \theta\in {\mathcal{T}}_{j+1} \,.
\end{equation}
The co-boundary operator $\partial^{\rm o}$ is defined from each ${\mathcal{T}}_j^{*}$ 
and takes its values in the space ${\mathcal{T}}_{j+1}^{*}$. The boundary operator makes 
decreasing the dimension of the chains while the co-boundary operator makes it increasing.

\noindent 
The co-boundary operator is a good tool to express the second Kirchhoff's law. 
We re-express the fundamental property (\ref{VDI}) in terms of the co-boundary operator:
\begin{equation}\label{DVI}
\left\langle \partial^{\rm o} V | i\right\rangle=0,~\forall V \in{\mathcal{T}}_0^{*},
~\forall i\in{\mathcal{T}}^1 \,.
\end{equation}
For each edge $a$ we introduce the potential differences in term of the potential values 
$V_s$ for each vertex $s$ and the incidence matrix $B$ as introduced in (\ref{Bord-Arete}):
\begin{equation}\label{DDP}
U_a=\sum_{s\in{\mathcal{T}}^0} B^s_a \, V_s \,.
\end{equation}
So we have: $ \partial^{\rm o}  V = \sum_{a \in {\mathcal{T}}^1} U_a
\left\langle a \right|$.
We introduce a closed circuit $\gamma$. We test the relation (\ref{DVI}) for 
$i = i_0 \sum_{a \in \gamma} \left| a \right\rangle $. Then 
$ \left\langle \partial^{\rm o} V | i\right\rangle> = i_0 \sum_{a \in {\cal T}^1}  \sum_{b \in \gamma}  
U_a \left\langle a | b  \right\rangle $ $ = i_0  \sum_{a \in {\cal T}^1}   U_a = 0   $. 
The mesh Kirchhoff's law express that the sum of the potential differences along a 
closed circuit is identically equal to zero.
For example if we take a look on circuit 2-3-4 figure 1, we have: $U_2-U_4-U_3=0$.

\bigskip \bigskip  \noindent  {\bf \large 7) \quad  Spanning tree for pair of nodes currents}

\noindent 
We assume now that the network ${\mathcal{T}}$ is connected. 
To fix the ideas we suppose more precisely that this network is simply 
connected {\it i.e.} does not contain any hole. When this hypothesis is not 
satisfied (a torus to fix the ideas) we refer to the contribution of 
Rapetti {\it et al.}  \cite{RDB03}.
A spanning tree $\mathcal{A}$ is a subgraph of the set of edges, doesn't 
contain any cycle, and is composed with a number of edges equal to the number 
of vertices minus one; if we add to this spanning tree an edge $a$ which doesn't 
belongs to $\mathcal{A}$ we obtain a cycle $\gamma$ composed by $a$ and edges of 
the spanning tree. We refer for a precise definition to the book of Berge  \cite{BERGE}.

\noindent 
Once a spanning tree ${\mathcal{A}}$ is fixed, any 1-chain can be decomposed 
in terms of boundary of faces plus a term associated to the spanning tree. 
In particular, each current can be decomposed in the previous form:
\begin{equation}\label{EC}
 \left| i \right\rangle = \sum_{f \in {\cal T}^2} \beta_f  \partial \left|  f \right\rangle +
\sum_{\alpha \in {\cal A}} \theta_{\alpha} \left| \alpha \right\rangle \,.
\end{equation}
The first term $\sum_{f \in {\cal T}^2} \beta_f  \partial \left|  f \right\rangle $ 
corresponds to the meshes currents in Kron's terminology, and the second one 
$\sum_{\alpha \in {\cal A}} \theta_{\alpha} \left| \alpha \right\rangle$ to the nodes-pairs 
currents. The formula (\ref{EC}) describes the direct sum of these two spaces. 
It is denoted as the ``complete space'' in Kron's approach.

\noindent 
We have the following theorem: if the current $i$ satisfies the Kirchhoff's law 
(\ref{VDI}) then the nodes pair currents is reduced to zero. We have in 
relation (\ref{EC}): $\theta_\alpha  = 0$ for each edge $\alpha$ of the spanning 
tree $\mathcal{A}$. The proof can be conducted as follows. Consider an arbitrary 
edge $\alpha \in {\mathcal{A}}$. We construct a potential $V$ as the one explicited 
on the figure 1.

\bigskip   
\centerline    { \includegraphics[width=.29 \textwidth]  {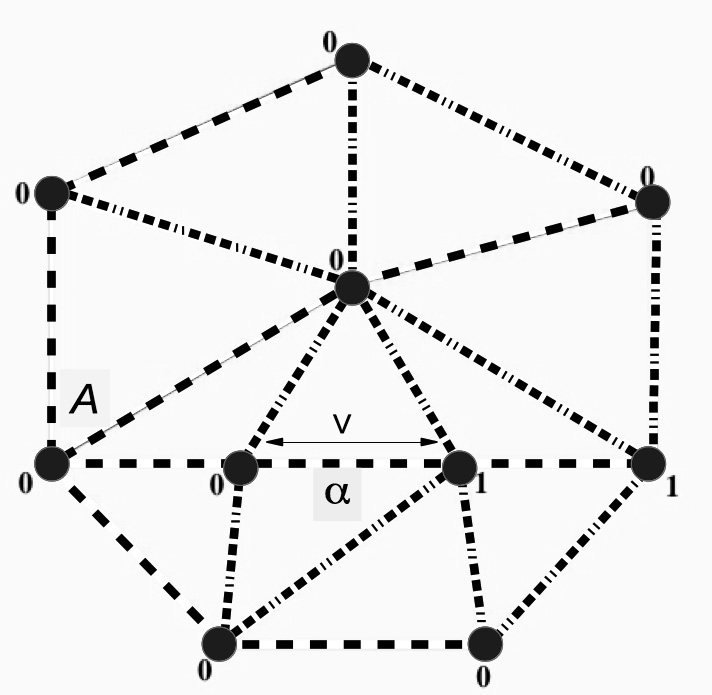} }  

\smallskip   \centerline  { {\bf Figure 1}. \quad Spanning tree. }
\bigskip   

\noindent 
We have $ \partial^{\rm o} V_\alpha = \left\langle \alpha \right|$ plus a sum 
related to edges that does not belong to the spanning tree $\mathcal{A}$.
Then we have the following calculus:  $0 = \left\langle V_\alpha |\partial 
i \right\rangle = \left\langle V_\alpha  |  \partial \big( \sum_{f \in {\mathcal T}^2} \beta_f 
\partial | f \rangle \big)\right\rangle $ $+ \left\langle V_\alpha  |  \partial \big(
   \sum_{\beta \in {\mathcal A}} \theta_{\beta} | \beta \rangle  \big)\right\rangle$ $
= \left\langle V_\alpha |  \partial \big(
   \sum_{\beta \in {\mathcal A}} \theta_{\beta} | \beta  \rangle  \big) \right\rangle$ $
= \left\langle  \partial^{\rm o}  V_\alpha |  \sum_{\beta \in {\mathcal A}} \theta_{\beta}  | 
\beta \right\rangle$  $= \left\langle \alpha | \sum_{\beta \in {\mathcal A}} \theta_{\beta}
  |  \beta \right\rangle$  $= \theta_{\alpha}$. The property is established.

\noindent 
When the Kirchhoff's law are satisfied, the electrical branches currents can be 
represented by the meshes currents $\left| i \right\rangle = \sum_{f \in {\cal T}^2} 
\beta_f  \partial \left|  f \right\rangle$. In the general case, when charges are 
injected to nodes, the Kirchhoff's laws (\ref{VDI}) are no more satisfied. The node 
pair currents $\sum_{\alpha \in {\cal A}} \theta_{\alpha} \left| \alpha \right\rangle$ 
is not equal to zero and represents these charges variations. In Maxwell's equations, 
the charge conservation has two terms: ``$div J$'' which is represented by the 
mesh currents and ``$\partial_t \rho$'' represented by the node pair currents.

\bigskip \bigskip  \newpage \noindent  {\bf \large 8) \quad  The fundamental space of faces }

\noindent 
In electrodynamic we want to make a relation between the meshes currents $i$ 
and some quantity coming from the flux $\Phi$. It translates the general 
relation of electrodynamic between currents and electromotive forces. This flux $\Phi$ 
is in relation with a magnetomotive force through $i$. The induced electrical current 
described in the previous section can be linked with the meshes currents through 
$\left| i \right\rangle = \sum_{f \in {\cal T}^2} \beta_f  \partial \left|  f \right\rangle$. 
The boundary operator being linear, we can write: $\left| i \right\rangle = 
\partial \left(\sum_{f \in {\cal T}^2} \beta_f \left|  f \right\rangle \right)$. 
This makes appear clearly the magnetic flux $\Phi$ given by:
\begin{equation}\label{phi}
 \Phi = \sum_{f \in {\cal T}^2} \beta_f \left|  f \right\rangle \,.
\end{equation}
The magnetic flux $\Phi \in {\mathcal{T}}^2$ is associated with the faces in the 
complex cellular ${\mathcal{T}}$. The meshes current $i$ being under this view a boundary current.

We define a dissipation operator $W$: $$ {\mathcal{T}}_2 \times {\mathcal{T}}_2 
\ni ~\left(\Phi,~\Phi'\right) \longmapsto W\left(\Phi,~\Phi'\right) \in \mathbb{R} .$$ 
It creates a positive defined quadratic form: $W\left(\Phi,~\Phi\right) \geq 0$ and 
$ W\left(\Phi,~\Phi\right)=0 \Rightarrow \Phi=0$.
This quadratic form $W$ generates a linear application: ${\mathcal{T}}_2 \ni  
\Phi \longmapsto \zeta \, \Phi \in {\mathcal{T}}_2^{*}$ so that:
\begin{equation}\label{fem1}
 \left\langle \zeta \, \Phi | \Psi \right\rangle = W\left(\Phi,~\Psi\right) 
~~\forall \Phi,~\Psi \in {\mathcal{T}}_2 \,.
\end{equation}
This linear application is nothing else than the impedance operator. On the faces 
base of ${\mathcal{T}}_2$ it creates the impedance matrix $Z$:
\begin{equation}\label{Z}
W \Big( {\sum_{f \in{\mathcal{T}}^2} \alpha_f \left| f \right\rangle ,~ \sum_{g \in{\mathcal{T}}^2} 
\beta_g \left| g\right\rangle} \Big)\equiv \sum_{f,g \in{\mathcal{T}}_2} Z_{fg}\alpha_f\beta_g \,.
\end{equation}
The energy source for the mesh space is given by the electromotive force $e$: for 
$\Phi \in {\mathcal{T}}_2,~e \in {\mathcal{T}}_2^{*}$ the dual product 
$\left\langle e | \Phi \right\rangle$ is well defined and points out $e$ as the dual 
source for faces as well as $\theta$ for the current source for nodes. The natural 
space for the electromotive force is cochain of degree~2.

Equilibrating sources and dissipation (in the general sense, {\it i.e.} used energy, 
losses or stored ones) we have: 
$$W\left(\Phi,~\Psi\right)=\left\langle e | \Psi\right\rangle,~~\forall \Psi \in {\mathcal{T}}_2$$ 
which means
\begin{equation}\label{fem2}
\zeta \, \Phi = e \,.
\end{equation}
Then the current can be obtained through: $$i=\partial \,\zeta^{-1}\, e$$ 
which is the topological expression for the Kron's  tensorial equation 
$$i^\mu = y^{\mu\nu} e_\nu \, .$$

\bigskip \bigskip  \noindent  {\bf \large 9) \quad  Resolution of networks in 
complete spaces }

\noindent \qquad \quad  {\bf \large using the Kron's method }

\noindent 
The method starts with a graph. This graph is a engineer view of a real system. 
Through homotopy, homology, surgery, the problem is projected on a graph \cite{HATCH}. 
In this operation, we start by finding a ST passing through the various remarkable points 
of the structure. This tree can be drawn on a sheet. Each point has its own connection 
with the original 3D space attached to the structure. Figure 2 shows such simple 
graph obtained from a filter.

\bigskip   
\centerline    { \includegraphics[width=.35 \textwidth]  {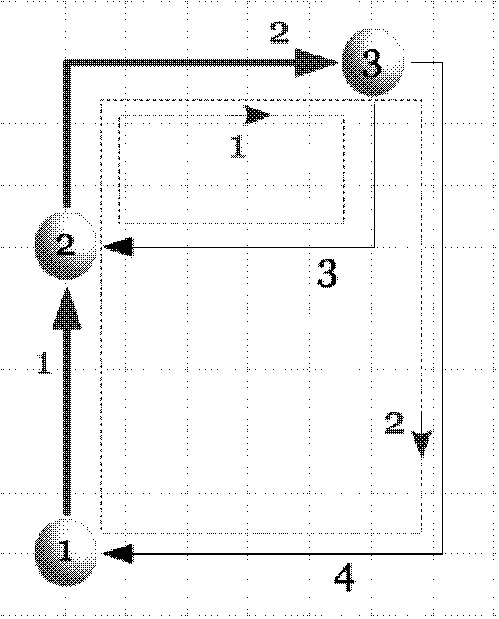}} 

\smallskip \centerline  { {\bf Figure 2}. \quad Graph of a filter. }
\bigskip   

\noindent 
On this graph, the ST is repaired by bold lines while closing edges are repaired by 
thin lines. Meshes are in doted and blue lines. Each edge current can be described 
depending on ST edges and meshes. For example, current of the edge 1 depends on ST 
edge 1 and mesh 2. The sign of this dependance is positive, as all currents flows 
in the same direction. When the meshes are constructed through closing edges from the 
ST, all the sign of the dependances are positive. This is a remarkable property linked 
with this construction method. The edge currents projection on ST edges and meshes is 
synthesized in a connectivity matrix C. If we choose to number the edges firstly as 
functions of closed edges belonging to the ST, then to meshes, the C matrix has a 
particular organization:
\begin{equation}
C=\left[
\begin{array}{cc}
Q & L \\ 0 & I
\end{array}
\right]
\end{equation}
As edges belonging to the ST are firstly numbered, they depend on both closing edges 
and meshes - that's why the submatrixes $Q$ and $L$ make links between these edges
 and both current sources and meshes. On the other hand, edges obtained by closing 
paths depend only on meshes. So a unity matrix links these edges with the meshes, 
and a zero matrix shows that there's no links between them and the closing paths.

\noindent 
Each edge $(a)$ has its own intrinsic property represented by an operator $z_{aa}$. 
This operator can depends on the current value on the edge (non linear one). 
This operator is a metric component which we talk about in the next paragraph. 
To present the Kron's method, we accept that this metric has the form:
\begin{equation}
z=\left[
\begin{array}{cc}
A & B \\ E & D
\end{array}
\right]
\end{equation}
We note $i^t$ the ST edges, $i^c$ the closing path edges, $J$ the current sources 
belonging to the ST, $k$ the mesh currents. The connectivity is:
\begin{equation}
\left[\begin{array}{c}
i^t\\i^c
\end{array}\right]=\left[\begin{array}{cc}
Q&L\\0&1
\end{array}\right]\left[\begin{array}{c}
J\\k
\end{array}\right]
\end{equation}
The Kirchhoff's law for any edges can be written \cite{PEIK} (it can be used 
for any physics):
\begin{equation}
\left[\begin{array}{c}
0\\S
\end{array}\right]=\left[\begin{array}{c}
V\\0
\end{array}\right] +z \left[\begin{array}{c}
i^t\\i^c
\end{array}\right]=\left[\begin{array}{c}
V\\0
\end{array}\right]+z\left[\begin{array}{cc}
Q&L\\0&1
\end{array}\right]\left[\begin{array}{c}
J\\k
\end{array}\right]
\end{equation}
S are the mesh sources on which we come back later, V the ST potential differences. 
By multiplying on the left by $C^T$ (index $T$ here is for transpose operation)  
we make appearing a bilinear transformation $C^T z C$. Noting $A',B',E',D'$ the 
components of this triple products (they are the component of the metric in the 
complete space: ST plus meshes), we obtain finally:
 \begin{equation}
\left[\begin{array}{c}
0\\S
\end{array}\right]=\left[\begin{array}{c}
Q V\\L V
\end{array}\right] +\left[\begin{array}{cc}
A' & B'\\ E' & D'
\end{array}\right]\left[\begin{array}{c}
J\\k
\end{array}\right]
\end{equation}
As previously demonstrated $L V=0$. The first equation resolved is: 
$S=E'J+D'k \Rightarrow k=\left(D'\right)^{-1}\left[S-E'J\right]$. Then, knowing k, 
the edges voltages of the ST current sources can be obtained: $Q V=-\left(A'J+B'k\right)$.

\bigskip \bigskip  \noindent  {\bf \large 10) \quad  About the metric } 

\noindent 
Usually, a metric is a matrix that can be half-positive, singular \cite{PEN}. 
In Kron's theory, as we want to manipulate complex operator through this tensor, 
we accept non strictly positive matrixes as metric and non symetric. Many mathematicians 
can object that we are finally far from a metric? But the notion stills very relevant 
as it describes a distance under preferential paths for currents in the topology, and 
it leads to the generalized power of the network studied. As said in \cite{PEK} page 309 
(this is a translation from French): \textit{A general theory for linear networks, 
symetrics or not, should be based on a metric space definition with fundamental tensor...}. 
The connectivity $C$ can be seen as a group transformation that belongs to SO 
( \cite{PEN} page 271).

\noindent 
We accept from now to call metric the fundamental tensor obtained from $\mu$. We mean 
by this that to create the Hodge's relation between the magnetic flux and the electromotive 
force (emf) we need to give ourselves a metric. This metric is like a rule to make a 
correspondance between the vectorial surface flux and the scalar emf. The metric makes 
the link between the first space of mesh currents and the dual quantity of emf obtained 
from an Hodge process.

\noindent 
From years, electronicians use the \textit{mutual inductance} to translate cross talk 
between both isolated circuits of a transformer. As we said before, emf and so 
self-inductances are for us deeply linked with the mesh space. Figure 3 shows the 
graph for a transformer cross talking two simple circuits.

\bigskip   
\centerline    { \includegraphics[width=.45 \textwidth]  {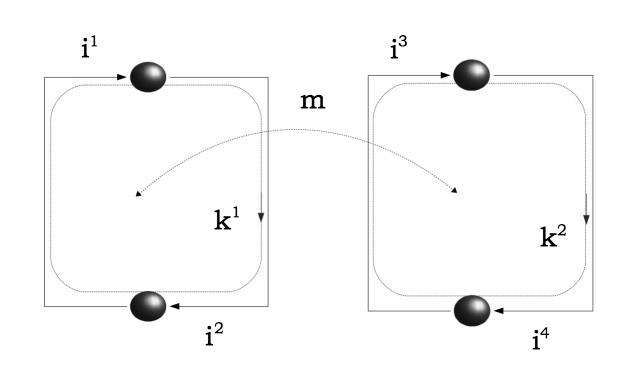}} 

\smallskip \centerline {  {\bf Figure 3}. \quad Graph of a transformer. }
\bigskip   

\noindent 
On the classical schematic presented fig. 3, we have two simple circuits having 
each of them a self-inductance and cross talked through a mutual inductance u. 
For each faces of each circuit, we can, following Maxwell's laws, compute the 
electric field circulation \cite{OMA}:
\begin{equation}
\oint_{\Gamma\equiv\partial f}\vec{E}\cdot\vec{d\Gamma}=\int_{A\in\Gamma}d\Gamma \frac{J}
{\sigma}+\int_{B\in\Gamma}d\Gamma \left[\int_t dt\frac{J}{\epsilon}\right]-hE_0
\end{equation}
A and B begin part of the boundary $\partial f$ of the face f, $hE_0$ begin the 
source electrical work giving energy to the circuit. Each kind of energy (dissipative 
through $\sigma$, potential through $\epsilon$) can belong to a single edge in the 
edge space. At the beginning of the circuit description, we have separate information 
for edges. Always for one circuit, one could be dedicate for dissipative energy 
(they are the classical resistors in electrical circuits) and another for potential 
energy (capacitances). Each elementary circuit can be for example a RC one, one edge 
A being a resistor, and another edge B being the capacitor. When we construct 
both circuits, we first begin to associate couple of edges in a single mesh. Using 
the mesh connectivity: $L^T=\left[\begin{array}{cc}1&1\end{array}\right]$, we obtain 
a first expression of the metric in the mesh space given for one circuit $g_n$ by:
\begin{equation}
g_n=L^TzL=\left[\begin{array}{cc}
1&1
\end{array}\right]\left[\begin{array}{cc}R_n & 0 \\ 0 & \frac{1}{C_np}
\end{array}\right]\left[\begin{array}{c}
1\\1
\end{array}\right]
\end{equation}
(p is the Laplace's operator, and $R_n$ and $C_n$ the resistance and capacitance 
of the circuit n). The metric for the two circuits comes from the union of all 
the metrics of separate circuits involved:  $g=\bigcup_n g_n$. For us it gives :
\begin{equation}
g=\left[\begin{array}{cc}R_1+\frac{1}{C_1p} & 0 \\ 0 & R_2+\frac{1}{C_2p}
\end{array}\right]
\end{equation}
As inductances belong only to the mesh space, we must add them to the previous 
matrix through a new one $\mu$, one including the inductance parts $L_n$:
\begin{equation}
g=g+\mu,~\mu=\left[\begin{array}{cc}L_1p & 0 \\ 0 & L_2p
\end{array}\right]
\end{equation}
Now if the first circuit has its own energy source $hE_0$ - we don't care from 
where it comes, the second circuit has no self energy generator. But a cross talk 
creates in its mesh an emf $e$. This emf comes from the mmf $F$ of the first circuit. 
Using the previous relations we can write finally a function between the mesh current 
of the first circuit $k^1$ and the emf of the second onde $e_2$: $e_2=-upk^1$. The 
cross talk is symetric, it means that the current in the second circuit $k^2$ creates 
an emf in the first one: $e_1=-upk^2$. Finally, the complete magnetic energy tensor 
added to the one obtain in the edges space becomes :
  \begin{equation}
\mu\left[\begin{array}{cc}L_1p & -up \\ -up & L_2p
\end{array}\right]
\end{equation}
We compute the emf added to the edge first circuit in a second step, after transformation 
of the edge metric in the mesh space and based on a function depending on the mesh 
current of another circuit. The mutual inductance appears here to be directly an 
application of the reluctance physic \cite{LeC}. A reluctance network and the 
associated graph can be construct where edges are tubes and vertex are meshes. 
This graph is the ``chord'' that is associated with the function u of mutual coupling. 
In next paragraph we will generalized this approach for Maxwell's fields. Before we 
obtain the extradiagonal component of our metric by another way. Another graph gives 
the same metric that the one of figure 3. Two spaces can be construct with these 
topologies. An isometric bijection exists between them that doesn't preserve the 
graph structure, but keep their common metric \cite{ENCu}. We consider the graph 
presented figure 4.

\bigskip   
\centerline    { \includegraphics[width=.45 \textwidth]  {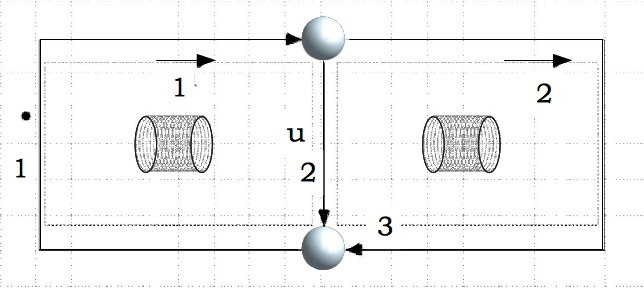}} 

\smallskip \centerline { {\bf Figure 4}. \quad Isometric graph. }
\bigskip   

\noindent 
Always without ST sources, the current can be expressed as: 
$I=\beta_1 \partial\left|f_1\right\rangle+\beta_2\partial\left|f_2\right\rangle$. 
The boundaries of the faces can be developed on the edges, depending on the directions 
chosen: $\partial\left|f_1\right\rangle=\left|a_1\right\rangle+\left|a_{12}\right\rangle,
~\partial\left|f_2\right\rangle=\left|a_2\right\rangle-\left|a_{12}\right\rangle$ 
where $a_1$ is edge 1 on the graph, $a_{12}$ edge 2 and $a_3$ edge 3. $f_1$ and $f_2$ 
are meshes 1 and 2 (doted blue lines fig.5). If $z_1$, $z_2$ and $z_3$ are the 
three components of the metric for the three edges 1, 2 and 3, the metric tensor 
in the mesh space ($L_m:z\rightarrow g$) becomes:
\begin{equation}
g=\left[\begin{array}{cc}
z_1+z_2 & -z_2 \\-z_2 & z_2+z_3
\end{array}\right]
\end{equation}
Comparing with previous one, we see that: $z_1=R_1+\left(C_1p\right)^{-1}+L_1p-up$, 
$z_2=R_2+\left(C_2p\right)^{-1}+L_2p-up$ and $z_2=up$. The natural energy distribution 
leads to the lagrangian operator obtained in the mesh space using the Kron's 
description of the networks \cite{GAB}. Graphs from figure 3 and 4 have not the 
same topology, but they do have the same metric, that's why we call \textit{isometric} 
the transformation from one to the other.

\bigskip \bigskip  \noindent  {\bf \large 11) \quad  Extension of chords to radiative fields}

\noindent 
All our previous discussion is based on graphs to describe a topology. When linked with 
real objects, these objects are reduced to edges through some algebraic topology operations 
and concepts like homotopy, homology, etc. Representation of reality is made through agencement 
of edges as pieces. The whole object modeled is a group of $R$ connex networks parts of a global 
graph. Each of these networks is a set of $B$ edges, joined by $N$ vertexes. The major characteristic 
of the object is its  number of meshes $M$ by respect of the equivalent Poincar\'e's law for complex 
cellular: $M=B-N+R$. In each connex network transmission of energy is perfectly controled because 
it goes from one vertex to another. By this way, the mathematical concepts can be applied to 
Maxwell's fields quite easily, under the hypothesis that the field behave like a bounded volume 
able to be represented by an edge. This is true for near and evanescent fields - for which 
macroscopic modelling are resistors, capacitances, inductances, reluctances, but this is not 
true for the free radiative field that cannot radiate in a bounded volume independently of the 
distance. We recall the fundamental difference between these two kinds of fields, and then we 
apply our previous method to these fields.

\bigskip   \noindent  \qquad {\bf \large a) \quad  Radiated versus evanescent fields} 

\noindent 
Basic demonstration for photon starts from the potential vector in the 
Coulomb's gauge \cite{FEYI}. Under this gauge, there is a formal separation of the 
transverse part of the field with the other components. Transverse part of the field 
leads to the photon concept of quantum mechanics. All the evanescent parts are the 
longitudinal components of the potential and scalar vectors. They can be modeled 
using inductances, capacitances, etc., that are properties (and components of the 
metric) of edges or meshes. Strictly speaking, evanescent modes of the field are 
virtual photon \cite{SaN}. The big difference for us is that the free field cannot 
be enclosed in a bounded edge (or a cycle of bounded edge, i.e. a mesh). It radiates 
in the infinity space. Another remarkable property of the far field is the radiation 
resistance. On the edge radiating, the property is increased of this radiation 
resistance. It can be shown that the radiation resistance is intrinsic to the radiated 
edge. When enclosed in a shielded room, the same edge has its metric modified by the 
metallic walls of the room. Reflecting the radiated field, the radiation resistance 
disapears due to the back induction coming from the walls. This process allows to 
desmontrate that the radiation resistance  is natively existing and is modified by 
the environment \cite{FEYCO} \cite{MAUT}.

\bigskip   \noindent  \qquad  {\bf \large b) \quad Emf and chords process for free radiated fields}

\noindent 
From years, engineers use antennas very simply: an input impedance gives the antenna 
equivalent circuit seen from the electronics. A radiation diagram describes the radiation 
of the antenna in all the free space. A gain, gives the relation between the free field 
and the power received. This very efficient modelling can be translated in our topology 
description. Basic principles are well explosed in  \cite{KM}. We present these principles 
before to translate them in topology.

\bigskip   \noindent  \qquad   {\bf \large c) \quad Antennas principles} 

\noindent 
Depending on the antennas gain, the receiving antenna can integrate the energy
 coming from the emitting antenna. The total radiated power $S_r$ in free space 
(whatever the distance r from the emitter) is:
\begin{equation}
S_r=\frac{P_t G_t}{4\pi r^2} e^{-\tau p}
\end{equation}
$P_t$ is the power delivered to the emitter, $G_t$ is its own antenna gain. The 
last term makes sure of the causality. The available power in reception is given 
by: $P_r=S_r A_r$. $A_r$ is the effective surface of reception of the antenna. The 
relation between gain and effective surface is:
\begin{equation}
G_x \left(\theta,\phi\right)=\frac{4\pi}{\lambda^2} A_x \left(\theta,\phi\right)
\end{equation}
Both gain or effective aperture (the other name for the effective surface of reception 
or emission) are functions of the 3D space. Using all the previous relations we find 
what is the Friss' equation:
\begin{equation}
\frac{P_r}{P_t}=A_r\left(\theta,\phi\right)\left[\frac{e^{-\tau p}}{\lambda^2 r^2}\right] 
A_t\left(\theta',\phi'\right)
\end{equation}
Now we just have to make the links between the power and the topology. For that, 
fundamental concepts previously stated should be used.

\bigskip   \noindent  \qquad  {\bf \large d) \quad  From power to topology } 

\noindent 

The total power of the radiated field is obtain through the radiative resistance on an edge. 
This component represents all the losses due to the radiation. As inductance and capacitance 
are linked with evanescent and lamellar fields, the radiated and transverse one which leaves 
the circuit (we consider no losses due to Joule effect in the wires) is linked with a resistance. 
For example we have: $P_t=R_{11}k^1$ and $P_r=\frac{\left(e_2\right)^2}{R_{22}}$. These relations 

give us the fundamental coupling impedance for far field:
\begin{equation}
z_{21}=\frac{e_2}{k^1}=\sqrt{R_{11}R_{22}A_r\left(\theta,\phi\right)
\left[G_{\theta,\phi,\theta',\phi'}\right]A_t\left(\theta',\phi'\right)}
\end{equation}
where $G_{\theta,\phi,\theta',\phi'}$ is a Green's kernel for the Friis formula. Once more we see 
that $e_2$ derive from an integration on a face that belongs to ${\mathcal{T}}^{2}$. This kind 
of interaction can be generalized, whatever the modes of the field propagated. Take a look to figure 5.

\bigskip   
\centerline    { \includegraphics[width=.45 \textwidth]  {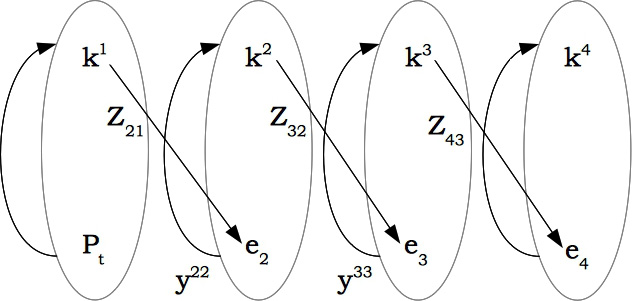}} 

\smallskip \centerline { {\bf Figure 5}. \quad General free field process. } 
\bigskip   

\noindent 
We see four simple meshes in interactions. For each of them we can identify the emf 
and mesh currents, and both emitted and received power. Going from the source of dynamic 
field $k^1$ to the induced emf $e_4$ we can describe the path: $e_4=z_{43}y^{33}z_{32}y^{22}z_{21}k^1$. 
The $z_{ij}$ are the coupling metric components and $y^{ij}$ the intrinsic inverse metric of 
meshes 2 and 3. $A_t$ and $A_r$ are properties attached to the graphs 1 and 4, i.e. to 
faces that belong to each graph of these circuits. In fact $e_2$ is the emf linked with $A_t$: 
it's a face resulting from a co-boundary applied to the cycle $k^1$. Each point on this 
face can be linked with a wave vector. The set of these wave vectors generate the flux $S_t$ 
(corresponding to $k^2$): the cotangent manifold of ${\mathbb{R}}^4$ generated by $k^1$ using 
an application $\phi$. The component of this set are transformed by the propagation operator 
$yzy$ to have its image through $S_r$ ($k^3$ and $\phi^{-1}$) where this time, $z$ is not a 
metric attached to the cellular complex ${\mathcal{T}}^\infty$ but to the 4D space-time fields
 propagation ${\mathbb{R}}^4$. The interface with the receiver is covered by the scalar 
product with $A_r\rightarrow k^3$. Then the emf $e_4$ is created by $k^3$.  Figure 6 shows 
the general process involved in the connection between the cell complex ${\mathcal{T}}^\infty$ 
and the 4D space ${\mathbb{R}}^4$

\centerline    { \includegraphics[width=.45 \textwidth]  {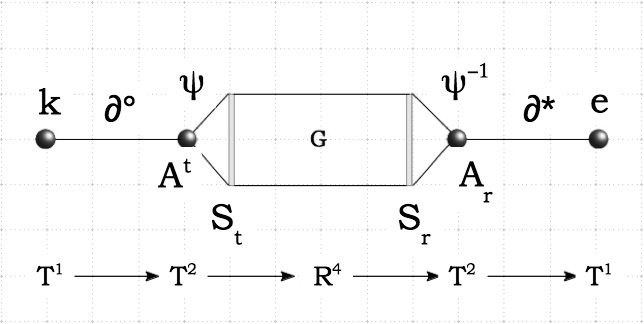}} 

\smallskip  \centerline  { {\bf Figure 6}. \quad Sequence of operation for free field.  }
\bigskip   

\bigskip \bigskip  \noindent  {\bf \large 12) \quad  Extension to functors} 

\noindent 
Previous function $\phi$ is a correspondance between both topological spaces 
${\mathcal{T}}^2$ and ${\mathbb{R}}^4$: $\phi:{\mathcal{T}}^2 \rightarrow {\mathbb{R}}^4$.  
We define as morphism of graphs, functions that preserve the graph structure, {\it i.e.}  keep the 
numbering of vertices, edges, etc. Cell complexes can be used to project any physical phenomenon 
on graphs. It was already usual to employ electroanalogy for various physics like thermal \cite{THE}, 
mechanics \cite{MECH}, quantum mechanics  \cite{SPROUL} and even biological information  \cite{BODY,MP}. 
In this last case, the author introduces the two categories of continuous and discrete spaces. 
This theory seems to be relevant for us, as it gives an algebraic approach for the method of 
chords  \cite{2loose}. A possible extension of his approach could be to generalize our previous 
work to functors between these two categories. By the way the method of chords replaces a large 
spatial domain where free electromagnetic energy propagates by a discrete link between topological 
objects. This kind of model could be extended to other physical configurations.

\bigskip \bigskip  \noindent  {\bf \large 13) \quad  A cavity problem as example}  

\noindent 
Considering a cavity with an aperture and an impeging parasitic wave, the problem 
is to compute the internal field induced by the external wave (figure 7). In the following, 
only the vertical polarization of the incident field will be considered.

\bigskip   
\centerline    { \includegraphics[width=.45 \textwidth]  {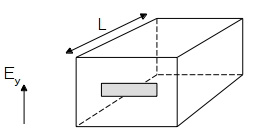}} 

\smallskip \centerline  { {\bf Figure 7}. \quad Electric field on cavity with an aperture.  } 
\bigskip   

\noindent 
This problem can be solved using the formalism developed above. Thanks to the analogy between 
respectively the electrical field and the voltage and then the magnetic field and the current, 
it is possible to define electrical equivalent scheme of the problem. In order to make such a 
representation, we will begin by the definition of different topological domain, corresponding 
to different regions that will be modeled separatly and then connected together to reconstruct 
the entire cavity. Three object have to be modeled:

$\bullet$  the incident field that can be assimilated to a voltage source in series with a 377 ohms impedance. 
This model represents the electrical field propagating in the free space that is equivalent to the 
mathematical representation in a plane wave shape of the incident field arriving with a normal incidence.

$\bullet$  the aperture, that will be seen from its middle point and that can be represented as two short 
circuited half transmission lines. The impedance seen from the center has already been extensively 
studied by different authors. The famous formulation given by Gupta (\cite{GUP}) has proved to 
efficiently represent the aperture impedance by a simple formula:
\begin{equation}
Z_{a}=120\pi^{2}[ln(2\frac{1+\sqrt[4]{1-(w_{e}/b)^{2}}}{1-\sqrt[4]{1-(w_{e}/b)^{2}}})]^{-1}
\end{equation}

$\bullet$  the cavity that can be seen as a transmission line along the main axis with some particular 
terminal conditions. On the side containing the aperture, the line will be directly connected 
in parallel with the aperture impedance model.The opposite side will be short circuited. In the 
present case, only the fundamental TE propagating mode will be considered. If we want to have 
results in a more large frequency bandwidth, we will have to consider each mode as an individual 
transmission line and to connect them in parallel.
In order to complete our model, we will add the possibility to realize a measurement in one point 
in the internal cavity. As in a real experiment, we will introduce a sensor that is modeled here 
as a resistor having a high value in order to avoid field perturbations inside  the cavity. This 
can be made by cutting the line into two half transmission lines located on each side of the 
transmission line. 

\noindent 
After having introduced the topology of our system, it is now important to give the Kron's 
transmission line model. In fact the more easier way to represent a transmission line is to 
give a quadrupole model. This quadrupole can be represented by two branches coupled by driven 
voltage sources. In fact, we can notice that such a representation is nothing else a circuit 
model of the impedance matrix of the line seen from its two extremities:
\begin{equation}
\left\lbrace\begin{array}{l}
V_1=Z_{11}i^1+Z_{12}i^2 \\ \\
V_2=Z_{21} i^1+Z_{22} i^2
\end{array}\right.
\end{equation}
For example, from the branch 1, the first equation shows two terms: the first one represents 
the influence of the current from this branch and the second one the effect of the current on 
the first branch. This coupling can be represented by a voltage source, the equation becomes: 
$V_{1}=Z_{11}I_{1}+V_{c}$. As a conclusion, the Kron's model of the cavity is given figure~8.

\bigskip   
\centerline    { \includegraphics[width=.45 \textwidth]  {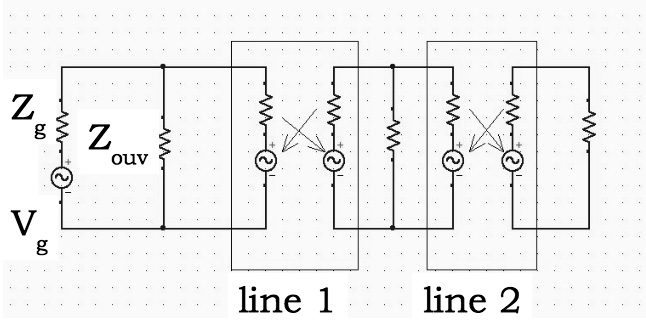}} 

\smallskip  \centerline  {  {\bf Figure 8}. \quad Kron's model for a cavity.   }
\bigskip   

\noindent 
In such an application, the important parameter that is usually modeled is the shielding 
effectiveness. This quantity was compared with the one obtained with a FDTD code. A very good 
agreement can be observed figure 9. Other examples between many others are given 
in  \cite{WCE2013},  \cite{ESA},  \cite{SMPS}.

\bigskip   
\centerline    { \includegraphics[width=.45 \textwidth]  {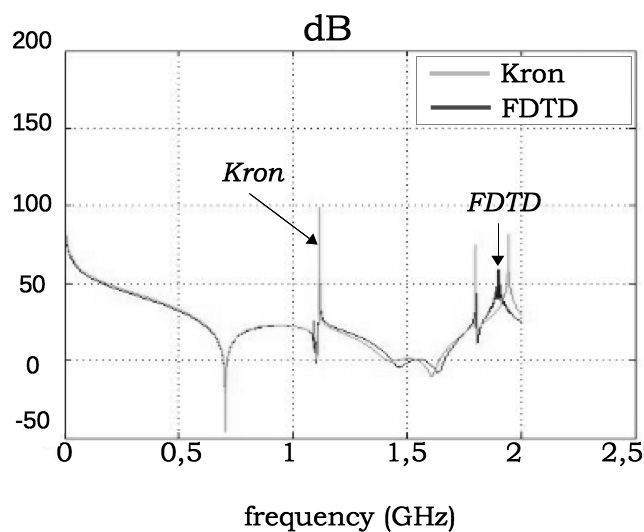}} 

\smallskip  \centerline  { {\bf Figure~9}. \quad Comparison with a  Finite Difference 
Time Domain computer code.  }
\bigskip   

\bigskip \bigskip  \noindent  {\bf \large 14) \quad  Antenna modelling}  

\noindent 
This second example shows the use of the formalism to compute the interaction between 
an antenna and a metallic wall. The objective of this experience was to understand the effect 
of a reflexion of energy on the radiation impedance of an antenna. One horn antenna is powered 
through an amplifier that delivers a amplitude modulated waveform at 10 GHz. Moreover, 50 cm in front 
of the antenna, there is a metalic wall or absorbers. The set-up of the experiments is shown 
figure~10.

\bigskip   
\centerline    { \includegraphics[width=.45 \textwidth]  {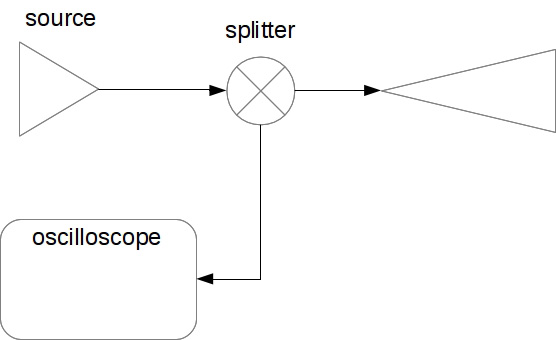}} 

\smallskip  \centerline  { {\bf Figure 10}. \quad Experiment set-up. }
\bigskip   

\noindent 
All the cable lengths were measured and the splitter resistors values characterized. 
Figure 11 shows the graph of the experience.

\bigskip   
\centerline    { \includegraphics[width=.45 \textwidth]  {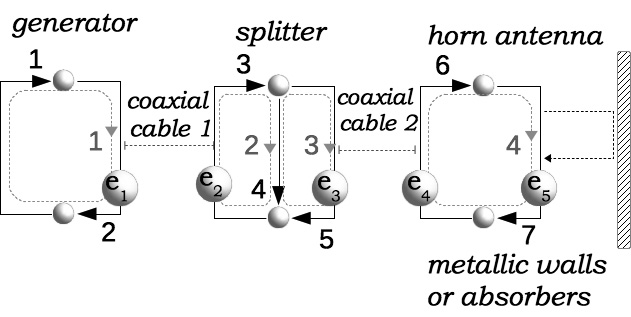}} 

\smallskip  \centerline  { {\bf Figure 11}. \quad Experiment equivalent graph.  }
\bigskip   

\noindent 
The splitter is made of three resistors of 17 ohms. Cables are simulated using Branin's model. 
For a cable of electrical length $\tau$, characteristic impedance $Z_c$, the Branin's model 
consists in two equations defining the electromotive force at each extremity of the line:
\begin{equation}
\left\lbrace \begin{array}{l}
e_1=\left(V_2-Z_c i^2\right) e^{-\tau p} \\ \\
e_2=\left(V_1+Z_c i^1\right) e^{-\tau p}
\end{array}\right.
\end{equation}
$V_1$ and $V_2$ are voltages respectively at the left and right of the line. $i^1$ and $i^2$ 
are the currents at the same extremities. Replacing $V_1$ and $V_2$ by their expressions 
depending on the loads and currents, Looking at the circuit figure 12, we obtain:
\begin{equation}
\left\lbrace \begin{array}{l}
V_1=E_0-R_0 i^1 \\ \\ V_2=RL i^2 \, . 
\end{array}\right.
\end{equation}

\bigskip   
\centerline    { \includegraphics[width=.45 \textwidth]  {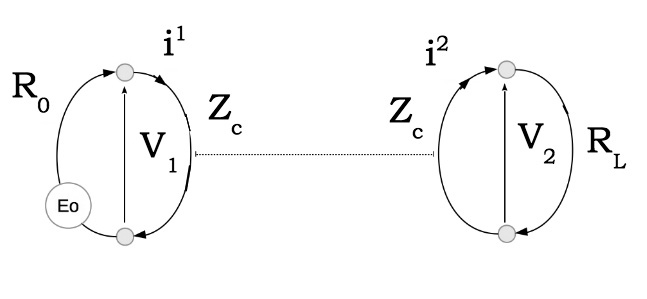}} 

\smallskip  \centerline  { {\bf Figure 12}. \quad Experiment equivalent graph.  }
\bigskip   

\noindent 
So, by replacement in (33) we understand that:

\begin{equation}
\left\lbrace \begin{array}{l}
e_1=\left(RL-Z_c \right) e^{-\tau p} i^2\\ \\
e_2-E_0e^{-\tau p} =\left(Z_c-R_0\right) e^{-\tau p} i^1
\end{array}\right.
\end{equation}
Any expression involving forms like $e_i/I^j$ can be replaced by an impedance interaction 
$z_{ij}$. Any line or guided wave structure can be replaced by an impedance tensor as:
\begin{equation}
\left[\begin{array}{cc}R_O+Z_c & \left(RL-Z_c \right)e^{-\tau p} \\ \\\left(Z_c-R_0\right)e^{-\tau p} & RL+Z_c
\end{array}\right]
\end{equation}

\noindent 
In the graph figure 11, two of these structures were used. The last edge of the graph 
represents the emitting horn antenna. It is a radiation resistance of fifty ohms. After a 
delayed time, when the field is reflected by the metallic wall in front of the antenna, a 
reflected field wave comes back in the horn and creates an electromotive force given by: 
$e_5=G\lambda \left(4\pi 2 R\right)^{-1}\sigma \sqrt{R_r} i^4$ ($R$ is the distance to the 
wall, $\sigma$ the reflection coefficient on the wall, $G$ the antenna gain and $R_r$ the 
radiation resistance). This creates an interaction given by $e_5/i^4$. The complete experience 
is detailed in \cite{OMHalG}.
Various measurements were made. Some with the wall equipped of a metallic plate. Some others 
with absorbers. The variation of radio frequency signal envelop at 10 GHz shows both cases of 
short circuit on the metallic wall or free space radiation on the absorbers. Difference between 
the signal compute under the Kron's formalism and measurements is of 1,2\%. This performance 
like others \cite{SUPant} was obtain using the method which allows to couple various accurate 
equations taken from different previous works.

\bigskip \bigskip  \noindent  {\bf \large 15) \quad  Conclusion }  

\noindent 
From fundamental definitions of discrete topology like the boundary operator and the 
notion of duality, we propose a mathematical model to formalize major results expressed by 
Gabriel Kron in his ``Tensorial Analysis of Networks''. One remarkable fact is that the 
chords introduced in a previous work appear as links between electromotive forces and meshes 
currents. We present an application of this result to an electromagnetic cavity without the 
help of three dimensional Maxwell's solver. Another application shows free radiated interaction 
between an antenna and a wall. This algebraic concept is flexible and the extension is 
under work for the modeling of multidisciplinary systems with networks.



\bigskip \bigskip \noindent {\bf \large  References } 

 


\begin{thebibliography}{99}


\bibitem{GKRONI}
G. Kron,  \textit{Tensorial Analysis of Networks}, General Electric editor, 1939.

\bibitem{GUP}
N.V. Balasubramanian, J.W.  Lynn, D.P. Sen Gupta,  
\textit{Differential Forms on Electromagnetic Networks}, Daniel Davey and Co., 1970.

\bibitem{K39}
G. Kron,  \textit{Non-Riemannian dynamics of rotating electrical machinery}, Romania,  1934.

\bibitem{K34}
A. Rottman,  \textit{Gabriel Kron et la formulation d'une technique de r\'esolution 
des syst\`emes complexes sur la base de la th\'eorie du circuit \'electrique}, 
Scientific journal AIM, 1988.

\bibitem{GKAST}
H.H. Happ,  \textit{Gabriel Kron and System Theory},  Union College Press, 1973.

\bibitem{BAH}
B. Hoffmann,  
``Kron's non Riemannian electrodynamics'', {\it Reviews of modern physics}, 
vol.~21, issue~3, p. 535-540, 1949.

\bibitem{ROTI}
J.P. Roth,  
``The validity of Kron's method of tearing'',
{\it Proceedings of the National Academy of Sciences of the United States of America}, 
vol.~41, issue 8, p. 599, 1955.

\bibitem{PEK}
M. Denis-Papin, A. Kaufmann,  \textit{Cours de calcul tensoriel}, Albin Michel, Paris, 1966.

\bibitem{HAPP}
H.H.  Happ,  \textit{Diakoptics And Networks}, Academic Press, 1971.

\bibitem{GK}
P.W. Gross, P.R. Kotiuga,  \textit{Electromagnetic  Theory and Computation}, MSRI editor, 2004.

\bibitem{BOS}
A. Bossavit,  \textit{Electromagn\'etisme en vue de la mod\'elisation}, Springer, New York, Berlin, 1993.

\bibitem{PEIK}
B. Peikari,  \textit{Fundamentals of network analysis and synthesis}, R.E. Krieger Pub. Co, 1982.

\bibitem{WHIT}
H. Whitney,  \textit{Geometric Integration Theory}, Dover publication, 2005.

\bibitem{RDB03}
F. Rapetti, F. Dubois,  A. Bossavit, 
``Discrete vector potentials for non-simply connected three-dimensional domains'', 
{\it SIAM journal on numerical analysis}, vol.~41, issue 4, p.~1505-1527, 2003.

\bibitem{BERGE}
C. Berge,  \textit{Theory of graphs and its applications}, Wiley, 1962.

\bibitem{HATCH}
A. Hatcher,  \textit{Algebraic Topology}, Cambridge university press, 2001.

\bibitem{PEN}
R. Penrose,  \textit{A la d\'ecouverte des lois de l'univers}, Odile Jacob - French translation, 2007.

\bibitem{OMA}
O. Maurice,  \textit{La compatibilit\'e \'electromagn\'etique des syst\`emes complexes}, 
Lavoisier, Paris, 2007.

\bibitem{LeC}
P. Lorrain, D.R. Corson,  \textit{Champs et ondes \'electromagn\'etiques}, 
Armand Collin, French version, page 426, 1970.

\bibitem{ENCu}
Encyclopedia universalis, \textit{Dictionnaire des math\'ematiques}, Albin Michel, page 651, 1997.

\bibitem{GAB}
R. Gabillard, \textit{Vibrations et ph\'enom\`enes de propagation}, Dunod, Paris, 1969.

\bibitem{FEYI}
R.P. Feynman,  \textit{Quantum Electro-Dynamics}, Addison-Wesley, 1961.

\bibitem{SaN}
A.A. Stahlhofen, G.  Nimtz,  \textit{Evanescent modes are virtual photons}, 
{\it Europhysics Letters}, vol.~76, issue 2, p.~189, 2006.

\bibitem{FEYCO}
R.B. Leighton, R.P.  Feynman, M.  Sands, M, \textit{Le cours de physique de Feynman}, 
InterEditions, 1979.

\bibitem{MAUT}
O. Maurice, \textit{Introduction d'une th\'eorie des jeux dans des topologies 
dynamiques}, Thesis, Limoges university, 2013.

\bibitem{KM}
J.D. Kraus, R.J. Marhefka,  \textit{Antennas}, Mc Graw Hill, 2002.

\bibitem{AReOM}
A. Reineix, O.  Maurice, P. Hoffmann, B. Pecqueux, P. Pouliguen, 
``Synthesis of the guided waves modelling principles under the tensorial analysis of network formalism'', 
European electromagnetics, EuroEM 2012.

\bibitem{THE}
J. Ouin, \textit{Transferts thermiques}, Casteilla edition, 1998.

\bibitem{MECH}
G.W. van Santen,  \textit{Vibrations m\'ecaniques}, Dunod edition, 1957.

\bibitem{BODY}
S. Grimnes, O.G. Martinsen,  \textit{Bioimpedance and bioelectricity basics}, Elsevier edition, 
second edition, 2008.

\bibitem{SPROUL}
R.L. Sproull,  \textit{El\'ements de physique moderne}, Masson editor, 1967.

\bibitem{MP}   
M. Poudret, J.P.  Comet, P. Le Gall, F.  K\'ep\`es, P.  Meseure, 
``Exploring Topological Modelling 
to Discriminate Models of Golgi Apparatus Dynamics'',  
European Conference on Complex Systems (ECCS), 2007.

\bibitem{2loose}
O. Maurice, A.  Reineix, P. Durand,  F. Dubois,  
``On mathematical definition of chords between networks'', European electromagnetics, EuroEM, 2012.

\bibitem{WCE2013}
P. Durand, O.  Maurice, A.  Reineix,  
``Generalized Interaction Principle Implemented in the Kron's Method'', 
Proceedings of the World Congress on Engineering, 2013.

\bibitem{ESA}
S. Leman, Y. Poir\'e, A. Reineix, F. Hoeppe, 
``Kron's method applied to the study of electromagnetic interference occuring in aerospace system'', 
Aerospace EMC, Proceedings ESA Workshop on EMC, 2012.

\bibitem{SMPS}
Y. Poir\'e, O. Maurice, M. Ramdani, M.  Drissi,  
``SMPS tools for EMI filter optimization'',  EMC Zurich2007,  18th International Zurich Symposium on EMC, 2007.

\bibitem{OMHalG}
O.  Maurice, 
``Exp\'erience  pour mettre en \'evidence la variation d'imp\'edance d'entr\'ee 
en fonction des interactions avec l'environnement'',  
HAL web docu\-ment\br 
 ``http://hal.archives-ouvertes.fr/docs/00/98/62/20/PDF/rapportSynthese\_manip\_ V2.pdf'', 2014.

\bibitem{SUPant}
A. Alaeldine, O. Maurice, J. Cordi,  R. Perdriau,  M. Ramdani, 
\textit{EMC-oriented analysis of electric near-field in high frequency},  ICONIC 2007.




\end{thebibliography}
\end{document}